\documentclass[12pt,oneside,english]{amsart}
\usepackage[latin1]{inputenc}
\usepackage{amssymb}

\makeatletter
\newtheorem{theorem}{Theorem}
\newtheorem{lemma}{Lemma}
\newtheorem{example}{Example}
\newtheorem{definition}{Definition}
\newtheorem{remark}{Remark}

\usepackage{babel}

\makeatother
\begin{document}
\baselineskip=17pt

\title[Exact exponent of  remainder term]{Exact exponent of  remainder term of Gelfond's digit theorem in binary case}

\author{Vladimir Shevelev}
\address{Department of Mathematics \\Ben-Gurion University of the
 Negev\\Beer-Sheva 84105, Israel. e-mail:shevelev@bgu.ac.il}


\begin{abstract}
We give a simple formula for the exact exponent in the remainder term of Gelfond's digit theorem in the binary case.
\end{abstract}

\maketitle

\section{Introduction}
 Denote for integer $m > 1$, $a\in[0, m-1]$.

\begin{equation}\label{1}
T^{(j)}_{m,a}(x)=\sum_{0\leq n< x,\enskip n\equiv a \mod m,\enskip s(n)\equiv j \mod 2}1,\qquad j=1,2
\end{equation}

where $s(n)$ is the number of 1's  in the binary expansion of $n$.

A.\enskip O.\enskip Gelfond \cite{5} proved that

\begin{equation}\label{2}
T_{m,a}^{(j)}(x) = \frac{x}{2m} + O(x^\lambda),\qquad j= 0,1,
\end{equation}

where

\begin{equation}\label{3}
\lambda= \frac{\ln 3}{\ln 4} = 0.79248125\ldots
\end{equation}

Recently, the author proved \cite{9} that the exponent $\lambda$ in the remainder term in (\ref{2}) is the best possible when $m$ is a multiple of 3 and is not the best possible otherwise.

In this paper we give a simple formula for the exact exponent in the remainder term of (\ref{2}) for an arbitrary $m$. Our method is based on constructing a recursion relation for the  Newman-like sum corresponding to (\ref{1})

\begin{equation}\label{4}
S_{m,a}(x) = \sum _{0\leq n < x,\; n\equiv a \mod m}(-1)^{s(n)},
\end{equation}
\newpage

It is sufficient for our purposes to deal with odd numbers  $m$ .
Indeed, it is easy to see that, if $m$ is even, then

\begin{equation}\label{5}
S_{m,a}(2x)= (-1)^a S_{\frac m 2,\lfloor\frac a 2\rfloor}(x).
\end{equation}

For an odd $m>1$, consider the number $r=r(m)$ of distinct cyclotomic cosets of 2 modulo $m$  [6, pp.104-105]. E.g., $r(15)=4$
since for $m=15$ we have the following 4 cyclotomic cosets of 2: $\{1,2,4,8\}, \{3,6,12,9\},\{5,10\},\newline\{7,14,13,11\}$.

Note that, if $C_1,\ldots, C_r$ are all different cyclotomic cosets of 2$\mod m$, then

\begin{equation}\label{6}
\bigcup^r_{j=1}C_j=\{1,2,\ldots, m-1\},\qquad C_{j_1}\cap C_{j_2}=\varnothing , \;\; j_1\neq j_2.
\end{equation}

Let $h$ be the least common  multiple of $|C_1|, \ldots, |C_r|$:

\begin{equation}\label{7}
h=[|C_1|,\ldots,|C_r|]
\end {equation}

Note that $h$ is of order 2 modulo $m$. (This follows easily, e.g., from Exercise 3, p. 104 in \cite{8}).

\begin{definition}  The exact exponent in the remainder term in (\ref{2})\enskip  is  $\alpha=\alpha(m)$ if
$$
T_{m,a}^{j}(x)=\frac{x}{2m}+O(x^{\alpha+\varepsilon}),
$$

and

$$
 T_{m,a}^{j}(x)=\frac{x}{2m}+\Omega(x^{\alpha-\varepsilon}),\qquad\forall\varepsilon>0.
$$

\end{definition}

Our main result is the following.

\begin{theorem}\label{t1}

If $m\geq 3$ is  odd, then the exact exponent in the remainder term in (\ref{2}) is

\begin{equation}\label{8}
\alpha=\max_{1\leq l \leq m-1}\left(1+\frac{1}{h\ln 2}\sum^{h-1}_{k=0}\left( \ln\left|\sin\frac{\pi l 2^k}{m}\right|\right) \right)
\end{equation}
\end{theorem}

Note that, if 2 is a primitive root of an odd prime $p$, then $r=1,\; h=p-1$. As a corollary of Theorem 1 we obtain the following result.
\newpage

\begin{theorem}\label{t2}
If $p$ is an odd prime, for which 2 is a primitive root, then the exact exponent in the remainder term in (\ref{2}) is

\begin{equation}\label{9}
\alpha=\frac{\ln p}{(p-1) \ln 2}.
\end{equation}
\end{theorem}

Theorem 2  generalizes the well-known result for $p=3$ (\cite{7}, \cite{2}, \cite{1}).
Furthermore, we say that 2 is a \slshape semiprimitive root \upshape modulo $p$  if 2  is  of order $\frac{p-1}{2}$ modulo $p$ and the congruence $2^x\equiv -1\mod p$ is not solvable. E.g., 2 is of order $8 \mod 17$, but the congruence
$2^x\equiv-1\mod 17$ has the solution $x=4$. Therefore, 2 is not a semiprimitive root of $17$. The first primes for which 2 is a semiprimitive root are (see[10], A 139035)

\begin{equation}\label{10}
7,23,47,71,79,103,167,191,199,239,263,\ldots
\end{equation}

For these primes we have $r=2,\;h=\frac{p-1}{2}$. As a second corollary of Theorem 1 we obtain the following result.

\begin{theorem}\label{t3}

If $p$ is an odd prime for which 2 is a semiprimitive root, then the exact exponent $\alpha$ in the remainder term in
(\ref{2}) is also given by (\ref{9}).
\end{theorem}

    In Section 2 we provide an explicit formula for $S_{m,a}(x)$, while in Sections 3-5 we prove  Theorems 1-3.

\section{Explicit formula for $S_{m,a}(x)$}

Let $\lfloor x \rfloor=N $. We have

$$
S_{m,a}(N)=\sum^{N-1}_{n=0,m|n-a}(-1)^{s(n)}=\frac 1 m \sum^{m-1}_{t=0}\sum^{N-1}_{n=0}(-1)^{s(n)} e^{2\pi i\frac{(n-a)t}{m}}
$$
\begin{equation}\label{11}
=\frac 1 m \sum^{m-1}_{t=0}\sum^{N-1}_{n=0}e^{2\pi i(\frac t m (n-a)+\frac 1 2 s(n))}.
\end{equation}

Note that the interior sum is of the form

\begin{equation}\label{12}
\Phi_{a, \beta}(N)=\sum^{N-1}_{n=0}e^{2\pi i(\beta (n-a)+\frac 1 2 s(n))},\qquad 0\leq \beta < 1.
\end{equation}

Putting

\begin{equation}\label{13}
F_\beta(N)=e^{2\pi i \beta a}\Phi_{a,\beta}(N),
\end{equation}
\newpage

we note that $F_\beta(N)$ does not depend on $a$.

\begin{lemma}\label{l1}
If $N= 2^{\nu_0}+ 2^{\nu_1}+ \ldots + 2 ^{\nu_\sigma},\;\;\nu_0 > \nu_r > \ldots >\nu_\sigma\geq 0$, then

\begin{equation}\label{14}
F_\beta(N)=\sum^\sigma_{g=0}e^{2\pi i(\beta \sum^{g-1}_{j=0} 2^{\nu_j}+\frac g 2)}\prod^{\nu_g-1}_{k=0}(1+e^{2\pi i (\beta 2^k +\frac 1 2)}).
\end{equation}
\end{lemma}

\bfseries Proof. \mdseries  Let $\sigma=0$. Then by (\ref{12}) and (\ref{13})

$$
F_\beta(N) = \sum^{N-1}_{n=0}(-1)^{s(n)}e^{2\pi i \beta n}
$$
\begin{equation}\label{15}
=1-\sum^{\nu_0-1}_{j=0}e^{2\pi i \beta 2^j}+\sum_{0\leq j_1 < j_2\leq\nu_0-1}e^{2\pi i \beta(2^{j_1}+2^{j_2})}-\ldots
\end{equation}
$$
=\prod^{\nu_0-1}_{k=0}(1-e^{2\pi i \beta 2^k}),
$$

which corresponds to (\ref{14}) for $\sigma=0$.

Assuming that (\ref{14}) is valid for every $N$ with $s(N)=\sigma+1$, let us consider $N_1=2^{\nu_\sigma}b + 2^{\nu_{\sigma+1}}$ where $b$ is odd, $s(b)=\sigma+1$ and $\nu_{\sigma+1}<\nu_\sigma$. Let

$$
N=2^{\nu_\sigma}b=2^{\nu_0}+\ldots +2^{\nu_\sigma};\;\; N_1=2^{\nu_0}+ \ldots +2^{\nu_\sigma}+2^{\nu_{\sigma+1}}.
$$

Notice that for $n\in[0, \nu_{\sigma+1})$ we have

$$
s(N+n)=s(N)+s(n).
$$

Therefore,

$$
F_\beta(N_1) =F_\beta(N)+\sum^{N_1-1}_{n=N}e^{2\pi i (\beta n+\frac 1 2 s(n))}
$$
\newline
$$
=F_\beta(N)+\sum^{\nu_{\sigma+1}-1}_{n=0}e^{2\pi i(\beta n+\beta N +\frac 1 2(s(N)+s(n)))}
$$
\newline
$$
=F_\beta(N)+e^{2\pi i (\beta N+\frac 1 2 s(N))}\sum^{\nu_{\sigma+1}-1}_{n=0}e^{2\pi i (\beta n+\frac 1 2 s(n))}.
$$

Thus, by (\ref{14}) and (\ref{15}),
\newpage

$$
F_\beta(N_1)=
$$
\newline
$$
\sum^\sigma_{g=0}e^{2\pi i(\beta\sum^{g-1}_{j=0}2^{\nu_j}+\frac g 2)}\prod^{\nu_g-1}_{k=0}(1+e^{2\pi i(\beta 2^k+\frac 1 2)}
$$
\newline
$$
+e^{2\pi i (\beta\sum^\sigma_{j=0}2^{\nu_j}+\frac{\sigma+1}{2})}\prod^{\nu_{g+1}-1}_{k=0}\left(1+e ^{2\pi i(\beta 2^k+\frac 1 2)} \right)
$$
\newline
$$
=\sum^{\sigma+1}_{g=0}e^{2\pi i (\beta\sum^{g-1}_{j=0}2^{\nu_j}+\frac g 2)}\prod^{\nu_g-1}_{k-0}\left(1+e^{2 \pi i (\beta 2^k +\frac 1 2)}\right).\blacksquare
$$

Formulas (\ref{11})-(\ref{14}) give an explicit expression for $S_m(N)$ as a linear combination of  products of the form

\begin{equation}\label{16}
\prod^{\nu_g-1}_{k=0}\left(1+e^{2\pi i (\beta 2^k+\frac 1 2)}\right),\;\; \beta=\frac t m,\;\; 0\leq t\leq m-1.
\end{equation}

\begin{remark}\label{r1}   One may derive  (\ref{14}) from a very complicated general formula of Gelfond \cite{5}. However, we prefered to give an independent proof.
\end{remark}

In particular,  if $N=2^\nu$,   then from (\ref{11})-(\ref{13}) and (\ref{15}) for

\begin{equation}\label{17}
\beta=\frac t m , \quad t=0,1,\ldots, m-1,
\end{equation}

we obtain the known formula cf. [3]:

\begin{equation}\label{18}
S_{m,a}(2^\nu)=\frac 1 m \sum^{m-1}_{t=1}e^{-2\pi i \frac t m a}\prod^{\nu-1}_{k=0}(1-e^{2\pi i\frac t m 2^k}).
\end{equation}

\section{Proof of Theorem 1}

     Consider the equation of order $r$

\begin{equation}\label{19}
z^r+c_1z^{r-1}+\ldots+c_r=0
\end{equation}

with the roots

\begin{equation}\label{20}
z_j=\prod_{t\in C_j}\left( 1-e^{2\pi i\frac t m}\right),\;\;\;j=1,2,\ldots,r.
\end{equation}
\newpage

Notice that for $t\in C_j$ we have

\begin{equation}\label{21}
\prod^{n+h}_{k=n+1}\left(1-e^{2\pi i\frac{t2^k}{m}}\right)=\left(\prod_{t\in C_j}\left(1-e^{2\pi i \frac t m}\right)\right)^{\frac{h}{h_j}}=z_j^{\frac{h}{h_j}},
\end{equation}

where $h$ is defined by (\ref{7}). Therefore, for every $t\in\{1,\ldots, m-1\}$, according to (\ref{19}) we have

$$
\prod^{n+rh}_{k=n+1}\left(1-e^{2\pi i\frac{t2^k}{m}} \right)
$$
\newline
\begin{equation}\label{22}
+c_1\prod^{n+(r-1)h}_{k=n+1}\left(1-e^{2\pi i\frac{t2^k}{m}} \right)+\cdots
\end{equation}
\newline
$$
+c_{r-1}\prod^{n+h}_{k=n+1}\left(1-e^{2\pi i\frac{t2^k}{m}} \right)+ c_r=0.
$$

After multiplication by $e^{-2\pi i \frac t m a} \prod^n_{k=0}\left(1-e^{2\pi i\frac{t2^k}{m}}\right)$ and summing over $t=1,2,\ldots, m-1$, by (\ref{18}) we find

\begin{equation}\label{23}
S_{m,a}\left(2^{n+rh+1}\right)+c_1 S_{m,a}\left(2^{n+(r-1)h+1} \right)+\cdots +c_{r-1}S_{m,a}\left(2^{n+h+1} \right)+c_rS_{m,a}\left(2^{n+1}\right)=0.
\end{equation}

Moreover, using the general formulas (\ref{11})-(\ref{14}) for a positive integer $u$, we obtain the equality

\begin{equation}\label{24}
S_{m,a}\left(2^{rh+1}u \right)+c_1S_{m,a}\left(2^{(r-1)h+1}u \right)+\cdots+c_{r-1}S_{m,a}\left(2^{h+1}u \right)+c_rS_{m,a}(2u)=0.
\end{equation}

Putting here

\begin{equation}\label{25}
S_{m,a}(2^u )=f_{m,a}(u),
\end{equation}

we have

\begin{equation}\label{26}
f_{m,a}(y+rh+1)+c_1f_{m,a}(y+(r-1)h+1)+\cdots +c_{r-1}f_{m,a}(y+h+1)+c_r f_{m,a}(y+1)=0,
\end{equation}
\newpage
where

\begin{equation}\label{27}
y=\log_2 u.
\end{equation}

The characteristic equation of (\ref{27}) is

\begin{equation}\label{28}
v^{rh}+c_1 v^{(r-1)h}+\cdots + c_{r-1} v^h +c_r=0.
\end{equation}

A comparison of (\ref{28})   and (\ref{20})-(\ref{21}) shows that the roots of (\ref{28}) are

\begin{equation}\label{29}
v_{j,w}=e^{\frac{2\pi i w}{h}}\prod_{t\in C_j}\left(1-e^{2\pi i \frac t m}\right)^{\frac 1 h},\;\; w=0,\ldots, h-1,\; j=1,2,\ldots, r.
\end{equation}

Thus,

\begin{equation}\label{30}
v=\max |v_{j,l}|= 2\max_{1\leq l\leq m-1}\left( \prod^{h-1}_{k=0}\left|\sin\frac{\pi l 2^k}{m}\right|\right)^\frac 1 h.
\end{equation}

Generally speaking, some numbers in (\ref{20}) could be equal. In view of (\ref{29}), the  $v_{j,w}$ 's have the same multiplicities.  If $\eta$ is the maximal multiplicity, then according to (\ref{27}), (\ref{30})

\begin{equation}\label{31}
S_{m,a}(u)=f_{m,a}(\log_2 u)= O\left( (\log_2 u)^{\eta-1}u^{\frac{\ln v}{\ln 2}}\right).
\end{equation}

     Nevertheless, at least

\begin{equation}\label{32}
S_{m,a}(u)=\Omega\left(u^\frac{\ln v}{\ln 2}\right).
\end{equation}

Indeed, let, say, $v=|v_{1,w}|$ and in the solution of (\ref{27}) with some natural initial conditions, all coefficients of $y^{j_1}v_{1,w}^y,\;\; j_1\leq \eta-1,\;\;w=0,\ldots, h-1$, are $0$. Then $f_{m,a}(y)$ satisfies a difference equation with the characteristic equation not having roots $v_{1,w}$ and the corresponding relation for

$S_{m,a}(2^n)$
(see (\ref{23})) has the characteristic equation (\ref{20}) without the root $z_1$. This is impossible since by (\ref{18}) and (\ref{21}) we have

$$
S_{m,a}(2^{h+1})=\frac 1 m \sum^r_{j=1}\sum_{t\in C_j}e^{-2\pi i\frac t m a}\prod^h_{k=1}(1-e^{2\pi i\frac t m 2^k})=
\frac 1 m \sum^r_{j=1}\sum_{t\in C_j}e^{-2\pi i\frac t m a} z_j^{\frac{h}{h_j}}.
$$
\newpage

Therefore, not all considered coefficients vanish, and (\ref{32}) follows. Now from (\ref{30})- (\ref{32})  we
obtain (\ref{8}).$\blacksquare$

\begin{remark}  In (\ref{8}) it is sufficient to let $l$ run over a system of distinct representatives of the cyclotomic cosets $C_1,\ldots,C_r$ of 2 modulo $m$.
\end{remark}

\begin{remark}
It is easy to see that there exists $l\geq 1$ such that $|C_l|=2$ if and only if $m$ is a multiple of 3. Moreover, in the capacity of $l$ we can take $\frac m 3$. Now from (\ref{8}) choosing $l=\frac m 3$ we obtain that $\alpha=\lambda=\frac{\ln 3}{\ln 4}$. This result was obtained in \cite{9} together with estimates of the constants in $S_{m,0}(x)=O(x^\lambda)$ and $S_{m,0}(x)=\Omega(x^\lambda)$ which are based on the proved in \cite{9} formula  $$S_{m,0}(x)=\frac 3 m S_{3,0}(x) + O(x^{\lambda_1})$$ for $\lambda_1=\lambda_1(m)< \lambda$  and Coquet's theorem \cite{2}.
\end{remark}

\begin{example}
Let $m=17,\; a=0$. Then $r=2,\; h=8$,

$$
C_1=\{1,2,4,8,16,15,13,9\},\;\; C_2=\{3,6,12,7,14,11,5,10\}.
$$

The calculation of $\alpha_l=1+\frac{1}{8\ln 2}\sum^{17}_{k=0}(\ln|\sin\frac{\pi l2^k}{17}|)$ for $l=1$ and $l=3$ gives

$\alpha_1=-0.12228749\ldots,\;\;\alpha_3=0.63322035\ldots$. Therefore by Theorem 1, $\alpha=0.63322035\ldots$. Moreover, we are able to prove that

$$
\alpha=\frac{\ln(17+4\sqrt{17})}{\ln 256}.
$$
\end{example}

Indeed, according to (\ref{23}), for $n=0$ and $n=1$ we obtain the  system $(S_{17,0}=S_{17})$:

\begin{equation}\label{33}
\begin{cases} c_1S_{17}(2^9)+c_2S_{17}(2)=-S_{17}(2^{17})\\c_1S_{17}(2^10)+c_2S_{17}(2^2)=-S_{17}(2^{18})
\end{cases}
\end{equation}

By  direct calculations we find
$$
S_{17}(2)=1,\;\; S_{17}(2^2)=1,\;\; S_{17}(2^9)=21,
$$
$$
S_{17}(2^{10})=29,\;\;S_{17}(2^{17})=697,\;\;\; S_{17}(2^{18})=969.
$$

Solving (\ref{33}) we obtain

$$
c_1= -34,\;\;\; c_2=17.
$$

Thus, by (\ref{23}) and (\ref{24})
\newpage

\begin{equation}\label{34}
S_{17}(2^{n+17})=34  S_{17}(2^{n+9})-17  S_{17}(2^{n+1}),\;\;\; n\geq 0,
\end{equation}

\begin{equation}\label{35}
S_{17}(2^{17}x)=34  S_{17}(2^9x)-17  S_{17}(2x),\;\;\; x\in \mathbb{N}.
\end{equation}

Putting furthermore

\begin{equation}\label{36}
S_{17}(2^x)=f(x),
\end{equation}

we have

$$
f(y+17)= 34  f  (y+9)- 17 (y+1),
$$

where $y=\log_2 x$. Hence,

$$
f(x)=O\left((17+4\sqrt{17})^{\frac x 8}\right),
$$

\begin{equation}\label{37}
S_{17}(x)=O\left((17+4\sqrt{17})^{\frac 1 8 \log_2 x}\right)=O(x^\alpha),
\end{equation}

where

$$
\alpha=\frac{\ln(17+4\sqrt{17})}{\ln 256}= 0.633220353\ldots
$$

\;\;\;\;\;\;\;\;\;\;

\section{Proofs of Theorems 2 and 3}

a) By the conditions of Theorem 2 we have $r=1,\;\; h=p-1$. Using (\ref{8}) we have

$$
\alpha=1+\frac{1}{(p-1)\ln 2}\ln\prod^{p-2}_{k=0}\left|\sin\frac{\pi 2^k}{p}\right|=1+\frac{1}{(p-1)\ln 2}\ln\prod^{p-1}_{l=1}
\sin\frac{\pi l}{p}.
$$

\;\;\;

Furthermore, using the identity [4, p.378],

$$
\prod^{p-1}_{l=1}\sin \frac{l \pi}{p}=\frac{p}{2^{p-1}}
$$

we find

$$
\alpha=1+\frac{1}{(p-1)\ln 2}\left(\ln p-(p-1)\ln 2 \right)=\frac{\ln p}{(p-1)\ln 2}.\blacksquare
$$
\newpage

\begin{remark}  In this case,  (\ref{24}) has the simple form
\end{remark}

$$
S_{p,a}(2^p u)+ c_1 S_{p,a} (2u) = 0.
$$

Since in the case of $a=0$  or 1 we have

$$
S_{p,a}(2)=(-1)^{s(a)},
$$

while in the case of $a\geq 2$,

$$
S_{p,a}(2 a) =(-1)^{s(a)},
$$

then putting

$$
u=\begin{cases} 1, \;\; a=0,1,\\
a,\;\; a\geq 2,\end{cases}
$$

we find
$$
c_1=(-1)^{s(a)+1}\begin{cases} S_{p,a}(2^p),\;\;  a=0,1, \\
S_{p,a}(a 2^p),\;\; a\geq 2 \end{cases}.
$$

In particular, if $p=3,\; a=2$ we have $c_1=S_{3,2}(16)=-3$ and

$$
S_{3,2}(8 u)= 3 S_{3,2}(2u).
$$

\begin{remark}
If  Artin's conjecture on the infinity of  primes for which 2 is  a primitive root  is true, then for $\alpha=\alpha(p)$ we have
\end{remark}

$$
\liminf_{p\rightarrow\infty}\alpha(p)=0.
$$

b) By the conditions of Theorem 3 we have $r=2,\;\;h=\frac{p-1}{2}$, such that for cyclotomic cosets of 2 modulo $p$

$$
C_1= -C_2.
$$

Therefore, in (\ref{8}) for $l_1=1$ and $l_2= p-1$ we obtain the same values. Thus,

$$
\alpha=1+\frac{2}{(p-1)\ln 2}\ln \left(\prod^{p-1}_{l=1}\sin\frac{\pi l}{p}\right)^\frac 1 2=\frac{\ln p}{(p-1)\ln 2}.\blacksquare
$$

\newpage

  Using Theorems 1-3, in particular we find $$\alpha(3)=0.7924..., \alpha(5)=0.5804..., \alpha(7)=0.4678...,\alpha(11)=0.3459,$$ $$\alpha(13)=0.3083...,\alpha(17)=0.6332...,\alpha(19)=0.2359...,
  \alpha(23)=0.2056...,$$ $$\alpha(29)=0.1734...,\alpha(31)=0.6358...,\alpha(37)=0.1447...,\alpha(41)=0.4339...,$$ $$\alpha(43)=0.6337...,\alpha(47)=
  0.1207... $$.

\end{document}